\documentclass[a4paper,11pt,twoside]{article}
\usepackage{amsmath,amssymb}

\newcommand{\eq}[1]{(\ref{#1})}
\newcommand{\dd}{\mathrm{d}} \newcommand{\e}{\mathrm{e}} \newcommand{\im}{\mathrm{i}}
\newcommand{\Tr}{\mathop{\mathrm{Tr}}}

\renewcommand{\Re}{\mathop{\mathrm{Re}}}
\newcommand{\bbbone}{{\mathchoice {\rm 1\mskip-4mu l} {\rm 1\mskip-4mu l}
{\rm 1\mskip-4.5mu l} {\rm 1\mskip-5mu l}}}

\newtheorem{theorem}{Theorem}

\newtheorem{prop}[theorem]{Proposition}

\begin{document}

\title{Stochastic differential equations for trace-class operators \\
and quantum continual measurements}
\author{A.~Barchielli and  A.~M.~Paganoni\\
{\small  {\sl Dipartimento di Matematica, Politecnico di Milano,}}\\
{\small{\sl Piazza Leonardo da Vinci 32, I-20133 Milano, Italy}} }
\date{}
\maketitle

\section{Introduction}\label{S1}

The theory of measurements continuous in time in quantum mechanics (quantum
continual measurements) has been formulated by using the notions of instrument
and positive operator valued measure \cite{Dav76}-\cite{Bar00}
%\cite{Dav76,BarLP83,Bel83,BarB91,BarP96a,Bar00},
arisen inside the \emph{operational approach} \cite{Dav76,Kra80} to quantum
mechanics, by using functional integrals \cite{BarLP82,Lup83,AlbKS97}, by using
quantum stochastic differential equations \cite{BarL85}-\cite{Bar97}
%\cite{BarL85,Bel88,Bel89a,BelS89,Bel89b,Bel90,Bel91a,Bel91b,Bel92a,Bel92b,BelS92,BarP96b,Bar97}
and by using classical stochastic differential equations (SDE's)
\cite{Bel88}-\cite{Bar00}.
%\cite{Bel88,Bel89a,BelS89,Bel89b,Bel90,Bel91a,Bel91b,BarB91,Bel92a,Bel92b,BelS92,%
%Bar93a,Bar93b,WisMil93,Bar94,BarH95,Kol95,Bar96,BarPZ98,BarZ98,Kol98,Doher99,BarP00,Bar00}.
Various types of SDE's are involved, and precisely linear and non linear
equations for vectors in Hilbert spaces and for trace-class operators. All such
equations contain either a diffusive part, or a jump one, or both.

In Section \ref{S2} we introduce a class of linear SDE's for trace-class
operators, relevant to the theory of continual measurements, and we recall how
such SDE's are related to instruments \cite{Dav76,Oza84} and \emph{master
equations} \cite{Lin76} and, so, to the general formulation of quantum
mechanics. In this paper we do not present the Hilbert space formulation of
such SDE's and we make some mathematical simplifications: no time dependence is
introduced into the coefficients and only bounded operators on the Hilbert
space of the quantum system are considered; for cases with time dependence see
for instance \cite{BarH95,BarPZ98,BarP00} and for examples involving unbounded
operators see for instance \cite{Bel89a}-\cite{BarZ98}.
%\cite{Bel89a,BelS89,Bel92a,BelS92,Kol95,Hol96,BarZ98}.
In Section \ref{S3} we
introduce the notion of \emph{a posteriori} state \cite{Bel83,Oza85} and the
non linear SDE satisfied by such states; then, we give conditions from which
such equation is assured to preserve pure states and to send any mixed state
into a pure one for large times. Finally in Section \ref{S4} we review the
known results about the existence and uniqueness of invariant measures in the
purely diffusive case \cite{Pag98,Kol98,BarP00} and we give some concrete
examples of physical systems. Other asymptotic results are given in
\cite{Kol95,Kol97}. In the whole presentation we try to underline the open
problems.

\section{Linear SDE's and instruments}
\label{S2}

Let us denote by $\mathfrak{L}(\mathcal{A}_1;\mathcal{A}_2)$ the space of
bounded linear operators from the Banach space $\mathcal{A}_1$ to the Banach
space $\mathcal{A}_2$ and let us set $\mathfrak{L}(\mathcal{A}_1)=
\mathfrak{L}(\mathcal{A}_1;\mathcal{A}_1)$. Let $\mathcal{H}$ be a separable
complex Hilbert space and let us consider the spaces of operators
$\mathfrak{L}(\mathcal{H})$, $\mathfrak{S}(\mathcal{H})=\big\{\rho \in
\mathfrak{L}(\mathcal{H}): \Tr\left\{\rho^* \rho\right\} < \infty \big\}$
(Hilbert-Schmidt operators), $\mathfrak{T}(\mathcal{H})=\big\{\rho \in
\mathfrak{L}(\mathcal{H}): \Tr\left\{\sqrt{\rho^* \rho}\right\} < \infty
\big\}$ (trace-class), with the norms $\|a\|_\infty$ (operator norm) for $a\in
\mathfrak{L}(\mathcal{H})$, $\|\rho\|_2= \sqrt{\Tr\{\rho^*\rho\}}$ for
$\rho\in\mathfrak{S}(\mathcal{H})$, $\|\rho\|_1= \Tr\{\sqrt{\rho^*\rho}\}$ for
$\rho\in\mathfrak{T}(\mathcal{H})$. We set $\langle a,\rho \rangle =
\Tr\left\{a^*\rho\right\}$, for $\rho \in \mathfrak{T}(\mathcal{H})$, $a \in
\mathfrak{L}(\mathcal{H})$ or for $\rho,\, a \in \mathfrak{S}(\mathcal{H})$.
Let us recall that $\mathfrak{L}(\mathcal{H})$ is the topological dual of
$\mathfrak{T}(\mathcal{H})$ and that $\|\cdot\|_\infty \leq \|\cdot\|_2\leq
\|\cdot\|_1$. Finally, \cite{Kra80} a quantum mechanical \emph{state} is a
trace-class operator such that $\rho = \rho^*$, $\rho \geq 0$, $\|\rho \|_1
\equiv \Tr\{\rho\}=1$; we denote by $\mathcal{S}(\mathcal{H})$ the set of all
states on $\mathcal{H}$. We denote by $\| \cdot\|_{1\to 1}$ the norm of a
bounded operator on $\mathfrak{T}(\mathcal{H})$ and we use similar notations
for operators on the other spaces.

Let $H$, $L_j$, $S_h$, $ j,h=1,2,\ldots$, be bounded operators on $\mathcal{H}$
such that $H=H^*$, $\sum_{j=1}^\infty L^*_jL_j=D_1$ and $\sum_{h=1}^\infty
S_h^*S_h=D_3$ are strongly convergent in $\mathfrak{L}(\mathcal{H})$. Let
$\mathcal{Y}$ be a locally compact Hausdorff space with a topology with a
countable basis, $\mathcal{B}(\mathcal{Y})$ be its Borel $\sigma$-algebra  and
$\nu$ be a finite measure on $\big(\mathcal{Y}, \mathcal{B}(\mathcal{Y})
\big)$. Let the operator  $\mathcal{J} \in \mathfrak{L}\big(
\mathfrak{T}(\mathcal{H}); L^1(\mathcal{Y},\nu;\mathfrak{T}(\mathcal{H}))\big)$
be such that the maps on $\mathfrak{T}(\mathcal{H})$ defined by
$\mathcal{R}_A[\rho] = \linebreak \int_A \mathcal{J}[\rho](y) \,\nu (\dd y)$,
$A\in \mathcal{B}(\mathcal{Y})$, are completely positive \cite{Lin76}; we set $
\mathcal{R}_\mathcal{Y}^{\,*}[\bbbone] = D_2 \in \mathfrak{L}(\mathcal{H})$.
According to \cite{BarP96a}, Def.\ 2 and Theor.\ 2, $\mathcal{R}_\cdot^{\,*}$
is a \emph{quasi-instrument}. Then, we introduce the following bounded
operators on $\mathfrak{T}(\mathcal{H})$:
\begin{equation} \label{2.1}
\begin{gathered}
\mathcal{L} = \sum_{i=0}^3 \mathcal{L}_i\, \qquad \mathcal{L}_0[\rho] = - \im
[H, \rho]\,, \qquad \mathcal{L}_1[\rho] =  \sum_{j=1}^\infty L_j \rho L_j^* -
\frac{1}{2} \bigl(\rho D_1 + D_1\rho\bigr),
\\
\mathcal{L}_2[\rho] = \mathcal{R}_\mathcal{Y}[\rho] - \frac{1}{2} \bigl(\rho
D_2 + D_2\rho\bigr), \qquad \mathcal{L}_3[\rho] =\sum_{h=1}^\infty S_h \rho
S_h^* - \frac{1}{2} \bigl(\rho D_3 + D_3\rho\bigr),
\end{gathered}
\end{equation}
The adjoint operators of $\mathcal{L}$, $\mathcal{L}_i$ are generators of
norm-continuous quantum dynamical semigroups \cite{Lin76}. We set also
\begin{equation}\label{2.3}
\mathcal{K}[\rho]=  \mathcal{L}_0[\rho]+  \mathcal{L}_1[\rho] - \frac 1 2 \,
D_2 \rho -\frac 1 2 \, \rho D_2  +\nu(\mathcal{Y})\rho + \mathcal{L}_3[\rho]\,.
\end{equation}

Let $\big( \Omega, (\mathcal{F}_t), \mathcal{F}, Q \big) $ be a stochastic
basis satisfying the usual hypotheses. Let $W_{j}(t)$ be continuous versions of
adapted, standard, independent Wiener processes with increments independent of
the past and $N (\dd y, \dd t)$  be an adapted Poisson point process on $
\mathcal{Y} \times \mathbb{R}_+$ of intensity $\nu(\dd y) \dd t$; $ N $ is
independent of the Wiener processes and with increments independent of the
past. We assume $(\Omega, \mathcal{F})$ to be a standard Borel space,
$(\mathcal{F}_t)$ to be the natural filtration of the processes $W,N$ and
$\mathcal{F}= \bigvee_{t\geq 0} \mathcal{F}_t$.

Let us now consider the following linear SDE for
$\mathfrak{T}(\mathcal{H})$-valued regular right continuous (RRC) processes:
\begin{equation}\label{2.2}\begin{split}
\sigma_t = \rho &+ \int_0^t \mathcal{K}[\sigma_{s^-}]\, \dd s +
\sum_{j=1}^\infty \int_0^t \bigl(L_j \sigma_{s^-} + \sigma_{s^-} L_j^* \bigr)
\dd W_j(s)
\\ {}&+
\int_{\mathcal{Y}\times (0,t]}\bigl(
\mathcal{J}[\sigma_{s^-}](y) - \sigma_{s^-} \bigr)N(\dd y, \dd s)\,,
\end{split}
\end{equation}
with nonrandom initial condition $\rho\in \mathfrak{T}(\mathcal{H})$; all the
integrals are defined in the weak, or $\sigma\big(\mathfrak{T}(\mathcal{H}),
\mathfrak{L}(\mathcal{H})\big)$, topology of $\mathfrak{T}(\mathcal{H})$.

The problem of the existence of a solution of eq.\ (\ref{2.2}) can be reduced
to a problem for SDE's in $\mathcal{H}$. It is possible, in many ways, to find
a larger probability space, where more Wiener and Poisson processes live, and
to construct a linear SDE for a process $\psi_t$ in $\mathcal{H}$ for which
existence and uniqueness of the solution can be proved by standard means and
such that the process $\sigma_t$, defined by $\langle a, \sigma_t \rangle =
\mathbb{E} [ \langle \psi_t |a \psi_t \rangle ]$ $\forall a \in
\mathfrak{L}(\mathcal{H})$, satisfies eq.\ (\ref{2.2}) \cite{BarH95,BarPZ98}.
From this representation one obtains that \emph{there exists a solution} of
eq.\ (\ref{2.2}) such that:
\begin{description}
\item[i.] the map $\rho \mapsto \sigma_t(\omega)$ is completely positive;
\item[ii.] if $\rho\in \mathcal{S}(\mathcal{H})$, then $\|\sigma_t\|_1=
\Tr\{\sigma_t\}$ is a positive martingale with $\mathbb{E}_Q\big[
\|\sigma_t\|_1 \big] = 1$.
\end{description}

About the uniqueness of the solution we are able to give some results under
additional assumptions (Theorem 1). It seems that simple estimates fail to give
uniqueness in the general case, so this remains an open problem. Let us stress
that even when existence and uniqueness can be proved directly from eq.
(\ref{2.2}), it seems difficult to prove the positivity property \textbf{i.}\
without going through the representation discussed above. Let us recall that
the solution must be an RRC process and that uniqueness is with respect to $Q$.

\begin{theorem}
(a) If $L_j=0$, $\forall j$, all the integrals in eq.\ (\ref{2.2}) are
meaningful in the norm topology of \ $\mathfrak{T} (\mathcal{H})$, the solution
is unique and also its existence follows directly from (\ref{2.2}).

(b) If $\sum_{j=1}^\infty \|L_j\|_\infty^{\,2} < +\infty$, \ $\sum_{h=1}^\infty
\|R_h\|_\infty^{\,2} < +\infty$, \ $\mathcal{J}[\rho](y) = \sum_{n=1}^\infty
J_n(y)\rho J_n(y)^* $, \ \ \ $J_n(y)\in \mathfrak{L}(\mathcal{H})$, \
$\int_\mathcal{Y} \Big( \sum_{n=1}^\infty \|J_n(y) \|_\infty^{\,2} \Big)^2
\nu(\dd y) <+\infty$, all the integrals in eq.\ (\ref{2.2}) are meaningful in
the norm topology of  $\mathfrak{S} (\mathcal{H})$, the solution is unique and
also its existence follows directly from (\ref{2.2}), seen as an equation in
the Hilbert space \ $\mathfrak{S} (\mathcal{H})$.
\end{theorem}
\textsl{Proof.} Case (a). All the statements follow, by standard arguments,
from the estimates:
\[
\Big\| \int_0^t \mathcal{K}[\sigma_{s^-}] \, \dd s \Big\|_1 \leq \|
\mathcal{K}\|_{1\to 1} \int_0^t \|\sigma_{s^-}\|_1\, \dd s\,,
\]
\[
\mathbb{E} \Big[ \Big\| \int_{\mathcal{Y}\times (0,t]} \bigl(
\mathcal{J}[\sigma_{s^-}](y) - \sigma_{s^-} \bigr) N(\dd y, \dd s)\Big\|_1
\Big] \leq \bigl(2 \|\mathcal{R_Y}\|_{1\to 1} + \nu(\mathcal{Y}) \bigr) \int
_0^t\mathbb{E}\big[ \|\sigma_{s^-}\|_1\big]\, \dd s\,.
\]

Case (b). The hypotheses on the operators $L_j$, $S_h$ give that $\mathcal{K}$
can be seen as a bounded linear operator on $\mathfrak{S}(\mathcal{H})$; so, we
have
\[
\mathbb{E}\Big[ \Big\| \int_0^t \mathcal{K}[ \sigma_{s^-}] \dd s \Big\|_2^{\,2}
\Big] \leq t \|\mathcal{K}\|_{2\to 2}^{\ 2} \int_0^t \mathbb{E} \big[ \|
\sigma_{s^-} \|_2^{\,2}\big] \dd s\,.
\]
The hypotheses on the $L_j$ give also
\[
\mathbb{E}\Big[ \Big\| \sum_j \int_0^t \big( L_j \sigma_{s^-} + \sigma_{s^-}
L_j^* \big) \dd W_j(s) \Big\|_2^{\,2} \Big] \leq 4 \sum_j \|L_j\|_\infty^{\,2}
\int_0^t \mathbb{E} \big[ \| \sigma_{s^-} \|_2^{\,2}\big] \dd s
\]
and the hypotheses on $\mathcal{J}$ give
\begin{multline*}
\mathbb{E}\Big[ \Big\|\int_{\mathcal{Y}\times (0,t]} \big(
\mathcal{J}[\sigma_{s^-}](y) - \sigma_{s^-} \big) N(\dd y, \dd s)
\Big\|_2^{\,2} \Big]
\\
\leq 2 \nu(\mathcal{Y}) t \Big( \int_{\mathcal{Y}}
\Big( \sum_n \|J_n(y)\|_\infty^{\,2}\Big)^2 \nu(\dd y) +\nu(\mathcal{Y}) \Big)
\int_0^t \mathbb{E} \big[ \| \sigma_{s^-} \|_2^{\,2}\big] \dd s\,.
\end{multline*}
Then, the uniqueness and the other properties follow by standard arguments.
\hfill{$\square$}

\smallskip

Let us denote by $\sigma_t^\rho$ the solution of eq.\ (\ref{2.2}) with initial
condition $\rho$ (if uniqueness does not hold, by $\sigma_t^\rho$ we mean
precisely that solution which has been constructed through a SDE on
$\mathcal{H}$). Then \cite{BarH95,BarPZ98} the equation
\begin{equation}\label{2.4}
\mathcal{I}_t(F)[\rho] = \mathbb{E}_Q[1_F\, \sigma_t^\rho ] \qquad \forall F\in
\mathcal{F}_t\,, \ \forall \rho \in \mathcal{S}(\mathcal{H})\,,
\end{equation}
defines a (completely positive) instrument \cite{Oza84} $\mathcal{I}_t$ with
value space $(\Omega,\mathcal{F}_t)$; the expectation of a trace-class operator
is defined in the weak topology of $\mathfrak{T}(\mathcal{H})$, as the
integrals in (\ref{2.2}). An instrument $\mathcal{I}$ is a measure such that
$\mathcal{I}(F)$ is a completely positive map on $\mathfrak{T}(\mathcal{H})$,
$\langle a, \mathcal{I}(\cdot)[\rho]\rangle $ is $\sigma$-additive and $\langle
\bbbone, \mathcal{I}(\Omega)[\rho] \rangle = \langle \bbbone, \rho \rangle$;
often it is the adjoint map of $\mathcal{I}$ which is called an instrument, as
in \cite{Oza84,BarH95,BarPZ98}.  Instruments give the most general setting for
representing measurement procedures in quantum mechanics: $\langle \bbbone,
\mathcal{I}(F)[\rho]\rangle $ is the probability of the event $F$ given the
premeasurement state $\rho$ and $\mathcal{I}(F)[\rho]\big/
\|\mathcal{I}(F)[\rho]\|_1$ is the postmeasurement state, conditional on the
event $F$. In our case we have a time dependent family of instruments, but by
property \textbf{ii.}\ the resulting probability measures are consistent and we
can define a unique probability $P_\rho$ on $(\Omega,\mathcal{F})$ by
\begin{equation}\label{2.5}
P_\rho(F) = \Tr \left\{ \mathcal{I}_t(F)[\rho] \right\} = \mathbb{E}_Q \left[
\left\|\sigma_t^\rho\right\|_1 1_F\right], \qquad \forall F\in \mathcal{F}_t\,.
\end{equation}

The interpretation of eqs.\ (\ref{2.4}) and (\ref{2.5}) is that
$\{\mathcal{I}_t,\, t\geq 0\}$ is the family of instruments describing the
continual measurement, the processes $W$, $N$ represent the output of this
measurement and $P_\rho$ is the physical probability law of the output.
Something more could be said on the instruments $\mathcal{I}_t$, linked to some
Markov property of $\sigma_t^\rho$ (see \cite{BarPZ98} Proposition 2.1).
Moreover, if we set
\begin{equation}\label{2.8}
m_j(t) = \langle  L_j + L_j^*,\rho_{t^-} \rangle, \qquad \lambda(t;y) = \langle
\bbbone, \mathcal{J}[\rho_{t^-}](y) \rangle ,
\end{equation}
it turns out \cite{BarH95,BarPZ98} that under the physical law $P_\rho$ the
processes
\begin{equation}\label{2.9}
\widetilde{W}_j(t) = W_j(t) - \int_0^t m_j(s) \dd s
\end{equation}
are independent standard Wiener processes and $N(\dd y,\dd t)$ is a counting
point process with stochastic intensity $\lambda(t;y)\nu(\dd y)\dd t$.

From eq.\ (\ref{2.5}) it follows that
\begin{equation}\label{2.6}
\eta_t = \mathcal{I}_t(\Omega)[\rho] =\mathbb{E}_Q [  \sigma_t^\rho]
\end{equation}
is the state to be attributed to the system at time $t$ if the output of the
measurement is not taken into account or not known; it can be called the
\emph{a priori} state at time $t$. It turns out that the \emph{a priori} states
satisfy the \emph{master equation}
\begin{equation}\label{2.7}
\frac{\dd\ }{\dd t}\, \eta_t = \mathcal{L}[\eta_t]\,, \qquad \eta_0=\rho\,.
\end{equation}
Master equations have been introduced in quantum mechanics as a way to
represent the dynamics of an open system, when ``memory" effects are
negligible.

\section{Nonlinear SDE's and \textit{a posteriori} states}
\label{S3}

If we introduce the random states
\begin{equation}\label{3.1}
\rho_t= \frac{\sigma_t^\rho}{ \left\|\sigma_t^\rho\right\|_1}\,,
\end{equation}
then we have, $\forall F\in \mathcal{F}_t$,
\begin{equation}\label{3.2}
\mathcal{I}_t(F)[\rho] =\mathbb{E}_Q[1_F \sigma_t^\rho]=
\mathbb{E}_{P_\rho}\left[1_F \frac{\sigma_t^\rho}{
\left\|\sigma_t^\rho\right\|_1}\right] = \mathbb{E}_{P_\rho}\left[1_F
\rho_t\right].
\end{equation}
According to \cite{Oza85}, $\rho_t(\omega)$ is a family of \emph{a posteriori}
states for the instrument $\mathcal{I}_t$ and the initial state $\rho$, i.e.\
$\rho_t(\omega)$ is the state to be attributed to the system at time $t$ when
the trajectory $\omega$ of the output is known, up to time $t$. Note that
$\eta_t = \mathbb{E}_Q [ \sigma_t^\rho ] = \mathbb{E}_{P_\rho}[\rho_t]$. It
turns out \cite{BarH95,BarPZ98} that the a posteriori states satisfy, under the
physical law $P_\rho$, the nonlinear SDE (again the integrals are defined in
the weak topology of the trace-class)
\begin{equation}\label{3.3}\begin{split}
\rho_t  = \rho &+ \int_0^t \mathcal{L}[\rho_{s^-}]\, \dd s + \sum_{j=1}^\infty
\int_0^t\bigl( L_j \rho_{s^-} + \rho_{s^-} L_j^* - m_j(s) \rho_{s^-}\bigr) \dd
\widetilde{W}_j(s)
\\  {}& +
\int_{\mathcal{Y}_t^\lambda} \Bigl(
\frac{1}{\lambda(s;y)}\,\mathcal{J}[\rho_{s^-}](y) - \rho_{s^-} \Bigr) \bigl(
N(\dd y, \dd s)- \lambda (s;y)\nu(\dd y)\dd s \bigr),
\end{split}
\end{equation}
where $\mathcal{Y}_t^\lambda= \{(y,s): y\in \mathcal{Y}, \, 0< s \leq t,\,
\lambda (s;y)> 0\}$. Let us stress that eq.\ (\ref{3.3}) has a solution by
construction, but that uniqueness is again an open problem, even under some
restrictive assumption; the point is that now the law of the noises depends on
the solution of the equation.

An interesting problem is to look for conditions ensuring the a posteriori
states $\rho_t$ for the instruments $\mathcal{I}_t$ to be almost surely (a.s.)
pure when the premeasurement state $\rho$ is pure; we recall that in the convex
set $\mathcal{S}(\mathcal{H})$ the pure states are the one-dimensional
projections. An interesting information-theoretical characterization of
instruments preserving pure states has been given in \cite{Oza86} and the
structure of this class of instruments has been found in \cite{Oza95}. A
measure of ``purity'' of a state $\rho$ is the so called linear entropy
$\Tr\{\rho (\bbbone -\rho)\}$; this quantity always belongs to the interval
$[0,1)$ and it is $0$ if and only if $\rho$ is a pure state. In our problem we
shall consider the linear entropy $g(t)= \langle\bbbone - \rho_t, \rho_t
\rangle$ and the mean linear entropy $G(t)= \mathbb{E}_{P_\rho}[g(t)]$ of the
solution of \eq{3.3}. For every random statistical operator $\rho_t$ we have $0
\leq G(t) < 1$; moreover, $G(t) = 0$ if and only if $\rho_t$ is a.s.\ a pure
state. So the study of the behaviour of the linear entropy is a way to analyze
whether the SDE \eq{3.3} preserves pure states or not. As in \cite{BarP00} we
can prove the following proposition.

\begin{prop}\label{puriinpuri}
Equation \eq{3.3} is ensured to preserve pure states, in the sense that
$\rho_t$ is a.s.\ a pure state for every a.s.\ pure initial condition, if and
only if $\mathcal{L}_3 = 0$ and there exist a set $A \subset \mathcal{Y}$, $A
\in \mathcal{B}(\mathcal{Y})$, a family of one-dimensional projections $P_y$,
$y\in A$, an operator $V \in \mathfrak{L}\big(\mathcal{H}; L^2(\mathcal{Y},\nu;
\mathcal{H})\big)$ such that
\begin{equation}\label{3.4}
\mathcal{J}[\rho](y) = \begin{cases} \sum_{\alpha}
\lambda_{\alpha}|(Vu_{\alpha})(y) \rangle \langle (Vu_{\alpha})(y)|\,, & y
\notin A\,,
\\[6pt]
\Tr\left\{\mathcal{J}[\rho](y)\right\}P_y\,, & y \in A\,,
\end{cases}
\end{equation}
where $\rho = \sum_{\alpha} \lambda_{\alpha} |u_{\alpha} \rangle \langle
u_{\alpha}|$ is a spectral decomposition of $\rho$.
\end{prop}
\textsl{Proof.} We observe that, by considering the embedding of
$\mathfrak{T}(\mathcal{H})$ into $\mathfrak{S}(\mathcal{H})$, the process
$\rho_t$ can be viewed as an Hilbert-space valued semimartingale (cf.\
\cite{Met82}, Def.\ 23.7). Then it is possible (cf.\ \cite{Met82}, Theor.\
27.2) to apply Ito's formula to the linear entropy $g(t)$ and to the mean
linear entropy $G(t)$ as done in Proposition 1 of \cite{BarP00}. In this way we
get the conditions $\mathcal{L}_3 = 0$ and $ \mathcal{J}[\rho](y)\big/
\Tr\left\{\mathcal{J}[\rho](y)\right\}$ is
$\bigl(\Tr\{\mathcal{J}[\rho](y)\}\nu(\dd y)\bigr)$--a.s.\ a pure state for
every pure state $\rho$. Then, by using Theorem 2 of \cite{BarP96a} and the
techniques of \cite{Oza95}, we get the representation (\ref{3.4}).
\hfill{$\square$}

\smallskip

Let us take now equation \eq{3.3} under the hypotheses of Proposition
\ref{puriinpuri} to guarantee that the equation preserve pure states; in
\cite{BarP00} we have also studied if it is possible to assure also that
\eq{3.3} map asymptotically mixed states into pure ones. Some examples of this
behaviour in the case of linear systems are given in \cite{Doher99}. We have
found some results only in the finite--dimensional case, the general case is
still an open problem. The proof of the next theorem is essentially as in
\cite{BarP00}.

\begin{theorem}\label{purificazione}
Let eq.\ \eq{3.3} preserve pure states and let $\mathcal{H}$ be
finite-dimensional. If it does not exist a bidimensional projection $P$ such
that
\begin{equation}\label{7.5}
\left\{
\begin{array}{ll}
P(L_j + L_j^*)P = z_j P \quad{} & \forall j
\\[6pt] P \mathcal{J}^*[\bbbone](y)P = q(y) P \quad{}& \nu\text{\rm -a.s.}
\\[6pt] P \mathcal{J}^*[\bbbone](y)P = 0 \quad{} &\nu|_A\text{\rm -a.s.}
\end{array}
\right.
\end{equation}
for some complex numbers $z_j$ and some complex function $q(y)$, then eq.\
\eq{3.3} maps asymptotically, for $t \to \infty$, mixed states into pure ones,
in the sense that for every initial condition the linear entropy vanishes for
long times:
\begin{equation}\label{entr}
\lim_{t \to \infty} \langle  \bbbone - \rho_{t},\rho_{t} \rangle = 0, \ \quad
{\rm a.s.}
\end{equation}
\end{theorem}
The space $\mathcal{H}$ being finite-dimensional, the vanishing of the linear
entropy is equivalent to the vanishing of the von Neumann entropy; so
(\ref{entr}) is equivalent to $\lim_{t \to \infty} - \Tr\{ \rho_t \ln \rho_t
\}$ ${}= 0$ a.s.\ \cite{Bar00}.

\section{Invariant measures}\label{S4}

Another interesting problem is the study of existence and uniqueness of an
invariant measure for SDE (\ref{3.3}). Here we treat the purely diffusive, pure
state preserving case, i.e. $J[\rho](y)=\rho$ (which implies $\mathcal{L}_2 =
0$) and $\mathcal{L}_3 = 0$; an example in the jump case is studied in
\cite{BarP00}. We assume also conditions (\ref{7.5}) to hold, in order to
guarantee that eq.\ (\ref{3.3}) map asymptotically mixed states into pure ones.
In this situation the support of any possible invariant measure is a subset of
the set of the pure states $M= \left\{\rho \in \mathfrak{T}(\mathcal{H}): \rho
= \rho^*, \rho^2 = \rho, \Tr{\rho} =1 \right\}.$ Moreover, we consider only the
case of a finite--dimensional Hilbert space; in such a case $M$ is a compact
subset of $ \mathfrak{T}(\mathcal{H})$. Since we work on a compact manifold, it
is useful \cite{IkW89} to use SDE's of Stratonovich type.

For dim$\,\mathcal{H} = n$, $\mathfrak{T}(\mathcal{H})$,
$\mathfrak{S}(\mathcal{H})$ and $\mathfrak{L}(\mathcal{H})$ coincide and reduce
to the space of $n\times n$ complex matrices. By taking a suitable orthonormal
basis, it is possible to represent every state as a $(n^2-1)$--dimensional real
vector and to represent $M$ as a $2(n-1)$--dimensional compact manifold
imbedded into $\mathbb{R}^{n^2-1}$ \cite{BarP00}. By working with such a
representation of states it is possible to obtain the Stratonovich form of eq.\
(\ref{3.3}); going back to the operator form we obtain
\begin{equation}\label{4.1}
\dd \rho_t = \mathcal{A}(\rho_t) \dd t + \sum_{j=1}^\infty
\mathcal{B}_j(\rho_t) \circ \dd W_j(t),
\end{equation}
\begin{equation}\label{4.2}\begin{split}
\mathcal{A}(\rho) &= - \im [H, \rho] + \sum_{j=1}^{\infty} \Big\{ \langle L_j +
L_j^*, \rho \rangle \mathcal{B}_j(\rho)+ {}
\\ {} &- \frac{1}{2} \big[(L_j + L_j^*) L_j \rho - \langle L_j +
L_j^*, L_j \rho \rangle \rho \big] + {} \\ {} & - \frac{1}{2} \big[ \rho
L_j^*(L_j + L_j^*) - \langle L_j + L_j^*, \rho L_j^* \rangle \rho \big]\Big\},
\end{split}
\end{equation}
\begin{equation}\label{4.3}
\mathcal{B}_j(\rho) = L_j \rho + \rho L_j^* - \langle L_j + L_j^*, \rho \rangle
\rho.
\end{equation}
Equation (\ref{4.1}) is a Stratonovich SDE on the $C_\infty$ manifold $M$; the
$\mathcal{A}, \mathcal{B}_j$ are $C_\infty$-vector fields. Moreover, it is
possible to show that, by some rotation on the Wiener processes, only a finite
number of Wiener processes really enter eq.\ (\ref{4.1}) (this is true because
$\mathcal{H}$ is finite-dimensional).

In \cite{Kol98} and \cite{BarP00} two theorems are given about existence and
uniqueness of an invariant measure; both the theorems are modifications of some
results in control theory \cite{Kli97}. For the notion of positive orbit of a
point for the deterministic control system associated to a diffusion see
\cite{Kli97} eq.\ (1.3).

\begin{theorem}\label{unicita'} \textrm{\rm\cite{BarP00}}
Let $M$ be a $d$-dimensional $C_\infty$ real compact manifold and let
$\mathcal{A}$, $\mathcal{B}_1, \ldots$, $\mathcal{B}_m$ be $C_\infty$ vector
fields; we consider the diffusion process defined by a Stra\-to\-no\-vich SDE
of type
\begin{equation}\label{elliptic}
\dd X_t = \mathcal{A}(X_t) \dd t + \sum_{j=1}^m \mathcal{B}_j(X_t) \circ \dd
W_j(t).
\end{equation}
If the generator $L$ of the diffusion process $X_t$ $\big(L = \mathcal{A} +
\frac{1}{2} \sum_{j=1}^m \mathcal{B}_j \mathcal{B}_j\big)$ is elliptic for
every $x \in M\backslash S$, where $S=\{x_1,\ldots,x_p\}$ is a finite set and
if in the points of $S$ the Lie algebra generated by the vectors fields
$\mathcal{A}, \mathcal{B}_1, \ldots, \mathcal{B}_m$ is full, i.e.\ on $S$
\begin{equation}\label{H}
\mathop{\mathrm{dim}} \mathcal{L}\mathcal{A}\{\mathcal{B}_0,\ldots,
\mathcal{B}_m\} = d\,,
\end{equation}
then there exists an unique invariant probability measure $\mu$ for equation
(\ref{elliptic}) and $\mathop{\mathrm{Supp}} \mu = M$.
\end{theorem}

\begin{theorem}\label{Kol} \textrm{\rm\cite{Kol98}}
Let $M$ be a $d$-dimensional $C_\infty$ real compact manifold and let
$\mathcal{B}_0, \mathcal{B}_1$ be $C_\infty$ vector fields; we consider the
diffusion process defined by a Stratonovich SDE of type (\ref{elliptic}) with
$m=1$. If there exists a smooth hypersurface $\Gamma$ in $M$ such that the
field $\mathcal{A}$ is transversal to $\Gamma$ and the equation $ \dot{X}_t=
\mathcal{B}_1(X_t)$ is a contraction outside $\Gamma$ (i.e. $\exists\, x_0 \in
M \setminus \Gamma$ such that $X_t$ tends to $x_0$ as $t \to \infty$ for all
initial points from $M \setminus \Gamma$), then, if we call $M_0$ the closure
of the positive orbit of $x_0$, all the points in $M_0$ are recurrent and all
the points in $M \setminus M_0$ are transient. Moreover there exists an unique
invariant probability measure $\mu$ and $\mathop{\mathrm{Supp}} \mu = M_0$.
\end{theorem}

Let us observe that Theorems \ref{unicita'} and \ref{Kol} deal with different
cases, in fact the number $m$ of fields entering the diffusive part of the SDE
is 1 for Theorem \ref{Kol}, while for Theorem \ref{unicita'} one necessarily
needs $m\geq d$.

In order to check the hypotheses of Theorem \ref{unicita'} we have firstly to
prove the ellipticity of the diffusion matrix on $M$ but a finite set of
points; the check of the ellipticity can be reduced to some properties of the
operators $L_j$. Firstly, the tangent space to $M$ in a point $\rho$ is $
T_{\rho}= \left\{ \tau \in \mathfrak{T}(\mathcal{H}): \tau = \tau^*,
\rho\tau+\tau\rho = \tau \right\}$; moreover for $\rho \in M$ it is possible to
write $\rho = |\psi \rangle \langle \psi |$ with $\psi \in \mathcal{H}, \
\|\psi \| =1$ and then it is possible to prove that we have also
$
T_{|\psi \rangle \langle \psi |}=\left\{ \tau \in \mathfrak{T}(\mathcal{H}):
\tau = |\psi\rangle \langle \phi| + |\phi \rangle \langle \psi |, \ \phi \in
\mathcal{H}, \ \langle \phi|\psi \rangle = 0 \right\}$. Secondly, by using
again the representation of the states in $\mathbb{R}^{n^2-1}$, one proves that
the diffusion matrix is elliptic in a point $\rho$ if and only if $\forall \tau
\in T_{\rho}$ there exists $j$ such that $\langle \tau , \mathcal{B}_j(\rho)
\rangle \neq 0$. Finally the ellipticity hypothesis of Theorem \ref{unicita'}
becomes: the diffusion matrix of the SDE (\ref{4.1}) is elliptic in a point
$|\psi\rangle \langle \psi|$ if and only if for all $\phi \in \mathcal{H}, \phi
\neq 0, \langle \phi|\psi \rangle = 0$ there exists $j$ such that
\begin{equation}\label{4.5}
\Re \langle \phi | L_j \psi \rangle \neq 0.
\end{equation}

For what concerns the application of Theorem \ref{Kol} to eq. (\ref{4.1}), one
has to study the equation $\dot{\rho_t} = \mathcal{B}_1 (\rho_t)$. If the
initial condition is $\rho_0 = |\psi_0 \rangle \langle \psi_0|$, $ \psi_0 \in
\mathcal{H}$, $\|\psi_0\| = 1$, then the solution of such an equation can be
written as
\begin{equation}\label{4.6}
\rho_t = \frac{|\psi_t\rangle \langle \psi_t|}{\|\psi_t\|^2}, \qquad \psi_t =
\e^{L_1 t} \psi_0.
\end{equation}
In \cite{Kol98} some concrete examples are studied with $L_1$ selfadjoint and,
in particular, with $L_1$ one-dimensional projection. Note that eq.\
(\ref{4.1}) is invariant for a change of sign of the $L_j$'s and the Wiener
processes; so, to check the hypotheses of Theorem \ref{Kol} also the equation
$\dot{\rho_t} = -\mathcal{B}_1 (\rho_t)$ can be considered.

Even for $\mathcal{H}=\mathbb{C}^2$, there are some physically interesting
examples. Let us consider a two--level atom of resonance frequency $\omega_0$,
stimulated by a monochromatic laser of frequency $\omega$; the interaction
between the atom and the electromagnetic field is mediated only by the
absorption/emission process. It is known that the detection schemes of the
fluorescence light known as homodyne/heterodyne detection and direct detection
can be treated by SDE's as eq.\ (\ref{3.3}), of diffusive and jump type
respectively \cite{WisMil93,Bar94,BarP96b}. In the case of an
homodyne/heterodyne detection scheme with a local oscillator in resonance with
the stimulating laser, the equation for the a posteriori states turns out to be
of the form (\ref{3.3}), with only the diffusive part and with
\begin{equation}\label{4.4}
L_j = \langle e_j | \alpha \rangle \sigma_-\,, \qquad H = - \frac{1}{2} \Delta
\omega \sigma_z + \im \langle \lambda| \alpha \rangle \sigma_- - \im \langle
\alpha| \lambda \rangle \sigma_+\,,
\end{equation}
where the sigma's are the Pauli matrices
\begin{equation*}
\sigma_+ =
\begin{pmatrix}
0&1\\ 0&0\\
\end{pmatrix},\qquad
\sigma_- = \begin{pmatrix} 0&0\\ 1&0\\
\end{pmatrix}, \qquad \sigma_z =
\begin{pmatrix}
1&0\\ 0&-1\\
\end{pmatrix},
\end{equation*}
$\Delta \omega = \omega -  \omega_0$, $\lambda, \alpha \in \mathcal{K}$, which
is a separable complex Hilbert space, $\{e_j, \ j=1,2,\ldots\}$ is a c.o.n.s.\
in $\mathcal{K}$; the choice of this c.o.n.s.\ depends on the concrete
implementation of the measurement scheme. The quantity $\|\alpha\|^2$
represents the natural line-width of the atom and $ 2 |\langle \alpha| \lambda
\rangle|$ is known as Rabi frequency; we assume both to be strictly positive,
i.e. $ \|\alpha\|^2 > 0$, $\Omega = 2 |\langle \alpha| \lambda \rangle |
> 0$.

In \cite{BarP00} we have proved that, if the complex numbers $ \langle e_j|
\alpha \rangle$ are not all proportional to a fixed one, then eq.\ (\ref{4.5})
is satisfied in the whole manifold $M$ but in the point $\rho_0 =
\begin{pmatrix}
0&0\\ 0&1\\
\end{pmatrix}$, where  the
condition (\ref{H}) is satisfied; so, in this case, the hypotheses of Theorem
\ref{unicita'} are fulfilled and there exists a unique invariant measure $\mu$
which has ${\text {\rm Supp}}\, \mu = M$.

In \cite{WisMil93} the problem of the invariant measure for a two-level atom is
studied by means of numerical simulations. The authors consider a heterodyne
detection scheme which satisfy all the hypotheses of Theorem \ref{unicita'}, as
we have shown in \cite{BarP00}, and a homodyne scheme which only in a
particular case can be handled via Theorem \ref{unicita'}; however, Theorem
\ref{Kol} can be applied. The model is characterized by $m=1$, $L_1= \e^{- \im
\phi}\, \|\alpha\| \sigma_-$, $\phi\in [0,2\pi)$, $\langle \alpha| \lambda
\rangle = \im \Omega / 2$. In this case the solution of the equation
$\dot{\rho_t} = \mathcal{B}_1 (\rho_t)$ can be easily computed via eq.\
(\ref{4.6}) and one obtains that the hypotheses of Theorem \ref{Kol} are
satisfied with $\Gamma = \emptyset$, $x_0=\rho_o=
\begin{pmatrix}
0&0\\ 0&1\\
\end{pmatrix}$. This means that there exists a unique invariant probability
measure $\mu$ with support $M_0$ implicitly given in Theorem \ref{Kol}. From
the structure of the fields $\mathcal{A}$ and $\mathcal{B}_1$, we conjecture
that the support of $\mu$ should be the whole $M$, apart from exceptional cases
such as the one studied in \cite{BarP00}, which corresponds to $\Delta \omega
=0$, $\phi=\pi/2$ and for which $M_0$ reduces to a circumference (for
$\mathcal{H}=\mathbb{C}^2$, $M$ can be represented as a spherical surface).

Let us end by discussing the connections between the invariant measure for eq.\
(\ref{3.3}) and the equilibrium state of the master equation (\ref{2.7}). The
following considerations apply both to the diffusive and to the jump case; we
assume that $\mathcal{H}$ is finite dimensional and that there exists a unique
invariant probability measure $\mu$ and that ${\text {\rm Supp}}\, \mu =:
M_0\subset M$. By definition of invariant measure we have
\begin{equation}\label{inv}
\int_M \mathbb{E}^\rho[f(\rho_t)] \mu(\dd \rho) = \int_{M_0} f(\rho) \mu(\dd
\rho)
\end{equation}
for every bounded measurable complex function $f$ on $M$; $\mathbb{E}^\rho$ is
the expectation in the case the initial condition for the process is $\rho$.
Moreover, by the uniqueness of $\mu$ we have ergodicity \cite{DapZab2}, i.e.
\begin{equation}\label{ergodic}
\lim_{T \to + \infty} \frac{1}{T} \int_0^T f(\rho_t) \dd t = \int_{M_0} f(\rho)
\mu(\dd \rho) \qquad \text{\rm a.s.}
\end{equation}
By applying these two facts to the function $f(\rho) = \langle a,\rho \rangle$,
where $a\in\mathfrak{L}(\mathcal{H})$, and by recalling that
$\mathbb{E}^\rho[\rho_t]=\eta_t$ satisfies the master equation (\ref{2.7}), we
obtain that the equation $\mathcal{L}[\eta] = 0$ has a unique solution solution
$\eta_{\rm{eq}}$ in $\mathcal{S}(\mathcal{H})$, given by
\begin{equation}\label{6.3}
\eta_{\rm{eq}} = \int_{M_0} \rho \mu(\dd \rho);
\end{equation}
moreover for every initial condition
\begin{equation}\label{6.4}
\lim_{T \to + \infty} \frac{1}{T} \int_0^T \rho_t \dd t = \eta_{\rm{eq}} \qquad
\text{\rm a.s.}
\end{equation}
Another quantity of interest is $D^2(a;\rho)= \langle (a^* - \langle
a,\rho\rangle) (a - \langle a^*,\rho\rangle),\rho\rangle \equiv \langle a^* a ,
\rho \rangle - |\langle a,\rho\rangle|^2$, $a\in \mathfrak{L}(\mathcal{H})$,
$\rho\in \mathcal{S}(\mathcal{H})$, which is a quantum analogue of the
variance. Note that we have the decomposition
\begin{equation}\label{6.5}
D^2(a;\eta_{\rm eq})= \int_{M_0}D^2(a;\rho)\mu(\dd \rho) + \int_{M_0} |\langle
a, \rho - \eta_{\rm eq} \rangle|^2 \mu(\dd \rho)
\end{equation}
and, by the ergodicity, we obtain
\begin{equation}\label{6.6}
\lim_{T\to +\infty} \frac 1 T \int_0^T D^2(a;\rho_t)\, \dd t=
\int_{M_0}D^2(a;\rho)\mu(\dd \rho) \,.
\end{equation}
For the case $m=1$, $L_1=L_1^*$, the asymptotic behaviour of $ D^2(L_1;\rho_t)$
has been studied in \cite{Kol98}.

\end{document}